\documentclass[10pt]{article}
\usepackage{makeidx}
\usepackage{multirow}
\usepackage{multicol}
\usepackage[dvipsnames,svgnames,table]{xcolor}
\usepackage{graphicx}
\usepackage{ulem}
\usepackage{amsmath}
\usepackage{amssymb}
\usepackage{float}
\usepackage{caption}
\author{Noah}
\title{}
\usepackage[paperwidth=612pt,paperheight=792pt,top=72pt,right=72pt,bottom=72pt,left=72pt]{geometry}

\makeatletter
	{\par\setlength{\parindent}{#3}
	\setlength{\leftmargin}{#1}       \setlength{\rightmargin}{#2}%
	\advance\linewidth -\leftmargin       \advance\linewidth -\rightmargin%
	\advance\@totalleftmargin\leftmargin  \@setpar{{\@@par}}%
	\parshape 1\@totalleftmargin \linewidth\ignorespaces}{\par}%
\makeatother

\newcommand{\qed}{\nobreak \ifvmode \relax \else
      \ifdim\lastskip<1.5em \hskip-\lastskip
      \hskip1.5em plus0em minus0.5em \fi \nobreak
      \vrule height0.75em width0.5em depth0.25em\fi}

\begin{document}

\title{Complex Objects in the Polytopes of the Linear State-Space Process}
\author{Noah E. Friedkin\textsuperscript{1}} 
\maketitle

\noindent \textit{ \textsuperscript{1} Center for Control, Dynamical Systems and Computation, College of Engineering, and Department of Sociology, College of Letters and Science,  University of California, Santa Barbara }.


\begin{abstract}
\noindent  A simple object (one point in $m$-dimensional space) is the resultant of the evolving matrix polynomial of walks in the irreducible aperiodic network structure of the first order DeGroot (weighted averaging) state-space process. This paper draws on a second order generalization the DeGroot model that allows complex object resultants, i.e, multiple points with distinct coordinates, in the convex hull of the initial state-space. It is shown that, holding network structure constant, a unique solution exists for the particular initial space that is a sufficient condition for the convergence of the process to a specified complex object. In addition, it is shown that, holding network structure constant, a solution exists for dampening values sufficient for the convergence of the process to a specified complex object. These dampening values, which modify the values of the walks in the network, control the system's outcomes, and any strongly connected typology is a sufficient condition of such control.               
\end{abstract}


\section{Introduction}

Finite convex sets are important in various areas of basic and applied mathematics, and appear in linear state-space processes. One area of application is the first-order DeGroot \cite{DeGroot1974} discrete-time state-space process in which the state of each point of a set $n$ points, ${\bf{x}}_{i}(k+1),\; i = 1,...,n$, is a convex combination of the immediately prior states of all points of the system ${\bf{x}}_{j}(k)  \in  \mathbb{R}^{n},\; j = 1,...,n$. The model was formulated as a  mechanism by which consensus might be reached among a set of individuals, and it has become the benchmark model of the literature on opinion dynamics. Its precursors include the models of French \cite{French1956} and Harary \cite{Harary1959}. In addition, the model has become increasingly prominent in control theory \cite{Olfati2004,Bullo2009,Huang2009,Moreau2005,Olfati2007}. When the mechanism unfolds in an aperiodic irreducible network, the system converges to a single value on the real number line and, more generally for an ${\bf{X}}(k) \in \mathbb{R}^{n \times m}$, to a single location in $m$-dimensional space. 

The second-order convex combination mechanism that is employed in this paper includes the DeGroot mechanism as a special case. It allows the emergence of complex objects composed of discrete points distributed in $m$-dimensional space. The emergent configurations of points include all polytopes (geometric objects with flat sides) and all other arrangements of points in the convex hull of the initial ${\bf{X}}(0)$ state of the system. In its application to the literature on opinion dynamics, these objects correspond to feasible patterns of influenced opinions that may include opinion clusters, each with a different mean opinion, and differentiated factions, each with a different consensus on an issue. More broadly, with respect to control theory applications, the model may be employed to move a set of $n$ points, in discrete steps, from initial positions in $m$-dimensional space to any specified configuration of points in the convex hull of ${\bf{X}}(0)$. The control occurs in the framework of a time-invariant network structure with an evolving matrix polynomial that is adjusted by a diagonal matrix of values. The diagonal values, specific to each point of the system, control the values of the walks in the network structure. In an irreducible structure, the detailed typology (its particular configuration of edges) does not constrain the feasible set of specified end-states of the system, nor do the magnitudes of initial distances among the points. 

\section{Second-Order Convex Combination Mechanism}      

Let ${\bf{X}}(0) \in \mathbb{R}^{n \times m}$ and ${\bf{V}}(k)  \in \mathbb{R}^{n \times n}$. Then ${\bf{X}}(k) = {\bf{V}}(k){\bf{X}}(0)$ is a set of points with $m$ coordinates in the convex hull of ${\bf{X}}(0)$ for all ${\bf{V}}(k)$ in the domain of nonnegative matrices with rows sums constrained to 1. Such ${\bf{V}}(k)$ are generated by the second order discrete-time state-space process \cite{Friedkin2006,FJ2011} 

\begin{equation}
{\bf{X}}(k+1) = {\bf{A}}{\bf{W}}{\bf{X}}(k) + ({\bf{I}}-{\bf{A}}){\bf{X}}(0), k=0,1,2,...,
\label{eq:FJ}
\end{equation}

\noindent where ${\bf{W}}$ is a nonnegative matrix with row sums constrained to 1, and ${\bf{A}}$ is a diagonal matrix constrained to ${\bf{0}} \le {\bf{A}} \le {\bf{I}}$. The process presents an evolving matrix polynomial 

\begin{align}
{\bf{V}}(k) & =  {\bf{AW}}{\bf{V}}(k-1)  + ({\bf{I}} - {\bf{A}}), \;\; {\bf{V}}(0) = {\bf{I}}, \;\; k=1,2,... , \\
& = ({\bf{AW}})^k + \big[{\bf{I}} + {\bf{AW}} + ({\bf{AW}})^2 + ... + ({\bf{AW}})^{k-1} \big] ({\bf{I}}-{\bf{A}}),\;k>2 \notag
\end{align}

\noindent which corresponds to walks in a network structure and preserves each ${\bf{V}}(k)$ as a nonnegative matrix with rows sums of 1: 

\begin{align}
{\bf{X}}(1) & = {\bf{V}}(1){\bf{X}}(0), \;\;{\bf{V}}(1) = {\bf{AW}} + {\bf{I}} -  {\bf{A}}, \notag \\
& = \big[{\bf{AW}} + {\bf{I}} -  {\bf{A}}\big]{\bf{X}}(0), \notag \\
{\bf{X}}(2) & = {\bf{V}}(2){\bf{X}}(0),\;\; {\bf{V}}(2) = {\bf{AW}}{\bf{V}}(1) + {\bf{I}} -  {\bf{A}}, \notag \\
& = \big[ {\bf{AW}}{\bf{V}}(1) + {\bf{I}} -  {\bf{A}} \big] {\bf{X}}(0), \notag \\
& =  \big[ ({\bf{AW}})^2 + \big({\bf{I}} + {\bf{AW}} \big)({\bf{I}}-{\bf{A}})  \big]  {\bf{X}}(0), \notag \\
& \vdots \notag \\
{\bf{X}}(k) & = {\bf{V}}(k){\bf{X}}(0),\;\; k > 2 \notag \\
& = \big[({\bf{AW}})^k + \big({\bf{I}} + {\bf{AW}} + ({\bf{AW}})^2 + ... + ({\bf{AW}})^{k-1} \big)({\bf{I}}-{\bf{A}}) \big]   {\bf{X}}(0), \notag \\
& \vdots \notag \\
{\bf{X}}(\infty) & = {\bf{V}}{\bf{X}}(0), \;\;  {\bf{V}}  \equiv \mathop{\lim }\limits_{k \to \infty } {\bf{V}}(k) \;\; 
\text{if this limit exists.}
\end{align}

\noindent Because $ {\bf{V}}(1) = {\bf{AW}} + {\bf{I}} -  {\bf{A}}$ is row-stochastic, all ${\bf{W}}{\bf{V}}(k-1)$ are row-stochastic, and all ${\bf{V}}(k) =  {\bf{AW}}{\bf{V}}(k-1)  + ({\bf{I}} - {\bf{A}}), \;\; k>1,$ are row-stochastic.  

The sequence  $\{ {\bf{V}}(k); k = 0,1,... \}$ converges if and only if the $\mathop{\lim }\limits_{k \to \infty } ({\bf{AW}})^k$ exists. The spectral radius $\rho({\bf{AW}})$ is controlled by  ${\bf{A}}$.  

\begin{itemize}
\item If ${\bf{A}} = {\bf{0}}$, then ${\bf{V}} = {\bf{I}}$,  and  ${\bf{X}}(\infty) = {\bf{X}}(0)$. 

\item If ${\bf{A}} = {\bf{I}}$, i.e, ${\bf{AW}}$ is stochastic,  then ${\mid {\lambda} \mid}_{\text{max}} = 1$ (Frobenius), the sequence  $ \{ {\bf{W}}^{k};k=0,1,,... \}$ converges, and $\{ {\bf{V}}(k);k=0,1,... \}$ converges, if and only if the eigenvalues of ${\bf{W}}$ for which $\mid {\lambda} \mid = 1$ are all $\lambda =  1$. 

\item If ${\bf{0}} \le {\bf{A}} \le {\bf{I}}$, i.e., ${\bf{AW}}$ is sub-stochastic, then ${\mid {\lambda} \mid}_{\text{max}} \le 1$, (Frobenius), the sequence  $ \{ ({\bf{AW}})^{k};k=0,1,,... \}$ converges, and $\{ {\bf{V}}(k);k=0,1,... \}$ converges (a) if the eigenvalues of ${\bf{AW}}$ for which $\mid {\lambda} \mid = 1$ are all $\lambda =  1$ or (b)  if  $\lambda =  1$ is not included in the spectrum. In the latter case, $( {{\bf{I}} - {\bf{AW}}})$ is nonsingular and  ${\bf{V}} = { ( {{\bf{I}} - {\bf{AW}}})^{ - 1}}({ {\bf{I}} - {\bf{A}} })$. 

\item If ${\bf{0}} \le {\bf{A}} < {\bf{I}}$, i.e, ${\bf{AW}}$ is strictly sub-stochastic, then ${\mid {\lambda} \mid}_{\text{max}} < 1$ (Frobenius), the sequence  $ \{ ({\bf{AW}})^{k};k=0,1,,... \}$ converges, and  $\{ {\bf{V}}(k);k=0,1,... \}$ converges to ${\bf{V}} = { ( {{\bf{I}} - {\bf{AW}}})^{ - 1}}({ {\bf{I}} - {\bf{A}} })$. In this case, $\mathop{\lim }\limits_{k \to \infty } ({\bf{AW}})^k = {\bf{0}}$, and ${\bf{V}} =  \big[ \sum_{k=0}^{\infty} ({\bf{AW}})^k \big]({\bf{I}} - {\bf{A}})$ involves the Newman series.     

\end{itemize}

In this framework, if ${\bf{W}}$ is aperiodic and irreducible, and ${\bf{A}} = {\bf{I}}$, then ${\bf{V}}  = \mathop{\lim }\limits_{k \to \infty } {\bf{W}}^{k} $. The Perron-Frobenius theorem applies and gives the \textit{simple} object, one point, with the coordinates ${\bf{X}}(\infty)$, in the convex hull of ${\bf{X}}(0)$. If ${\bf{W}}$ is aperiodic and irreducible, and ${\bf{0}} < {\bf{A}} < {\bf{I}}$, then the result  

\begin{equation}
{\bf{X}}(\infty)  =  \big( {\bf{I}} - {\bf{AW}} \big)^{-1} ({\bf{I}} - {\bf{A}}) {\bf{X}}(0)
\end{equation}
 
\noindent is a \textit{complex} object of discrete points, with the coordinates ${\bf{X}}(\infty)$. As ${\bf{A}}$ approaches ${\bf{I}}$, these complex objects converge to the simple object of the ${\bf{A}} = {\bf{I}}$ case. More generally, for a particular $a_{ii}$, $\mathop{\lim }\limits_{ {a_{ii}} \to 1 } = 1 $ and $\mathop{\lim }\limits_{ {a_{ii}} \to 0 } = 0 $. Hence there is no loss of generality, and considerable gain of parsimony, entailed in the restriction ${\bf{0}} < {\bf{A}} < {\bf{I}}$.   

Every ${\bf{X}}(0)$ is associated with a polytope that is a box of $m$-dimensions, defined by the maximum and minimum values of each column of ${\bf{X}}(0)$. All $i=1,..,n$ points of the system might be located at the vertices of the box.  The feasible initial coordinates of the $n$ points of the system include, as special cases, all convex polytopes in $m$ dimensions. These special cases occur when non-vacuous subsets of the $n$ points of the system occupy each vertex of a particular type of polytope (e.g., a pentagon) and all points are initially located at these vertices. Holding ${\bf{X}}(0)$ and ${\bf{W}}$ constant, the state-space process will generate complex objects  in the minimal convex set (convex hull) of ${\bf{X}}(0)$ for all $0 < {\bf{A}} < {\bf{I}}$.        

It is immediately apparent that: 

\begin{itemize}
\item for given $ \{{\bf{W}},{\bf{X}}(0) \}$, the diagonal values of ${\bf{A}}$ control the system's outcome; 
\item for given $\{ {\bf{W}}, {\bf{0}} < {\bf{A}} < {\bf{I}},{\bf{X}}(\infty) \}$, there exists a unique
\begin{align}
{\bf{X}}(0) & = {\bf{V}}^{-1}{\bf{X}}(\infty) \notag \\
 & =   ({\bf{I}} - {\bf{A}})^{-1}  \big( {\bf{I}} - {\bf{AW}} \big) {\bf{X}}(\infty) , 
 \label{eq:X0}
\end{align}
\noindent on which basis the state-space process will converge to the specified complex object ${\bf{X}}(\infty)$; and
\item for given $\{ {\bf{W}}, {\bf{X}}(0),{\bf{X}}(\infty) \}$, there may be no feasible ${\bf{0}} < {\bf{A}} < {\bf{I}}$ that satisfies  
\begin{equation}
{\bf{X}}(\infty) - {\bf{X}}(0) = {\bf{A}} \big[ {\bf{W}}{\bf{X}}(\infty) - {\bf{X}}(0) \big].
\end{equation}
\end{itemize}

An implication of equation \ref{eq:X0} is that, for a given ${\bf{W}}$, there exists an infinite number of  $\{{\bf{A}},{\bf{X}}(0)\}$ combinations with which the state-space process will converge to a specified complex object ${\bf{X}}(\infty)$. I draw on this implication in an analysis of the emergence of specified complex objects via the state-space process. I restrict the analysis to aperiodic irreducible ${\bf{W}}$. Such matrices correspond to strongly connected networks in which at least one path from $i$ to $j \ne i$ exists for all ordered $(i,j)$ pairs of nodes. Given a strong structure, a sufficient condition of aperiodicity is the existence of at least one positive value on its main diagonal, which corresponds to at least one positive resistance loop in the network. It is significant that the particular details of the typology of the network associated with ${\bf{W}}$, i.e., its idiosyncratic configuration of edges, do not matter. The detailed typology may be important in other respects, but it is not important with respect to reaching the end state of the process that is the specified complex object.     

The remainder of paper is organized into two main sections and concludes with a brief discussion:  

\begin{enumerate}
\item The first section illustrates various complex objects that may be formed under the structural conditions of an irreducible aperiodic ${\bf{W}}$. 
\item The second section addresses the particular combinations  $\{{\bf{A}},{\bf{X}}(0)\}$ with which the state-space process will converge to a specified complex object ${\bf{X}}(\infty)$ for a given ${\bf{W}}$. 
\end{enumerate}

\section{Complex Objects}
The complex objects that arise within the polytope associated with  ${\bf{X}}(0)$ are $n$ points with $m$ coordinates. If undirected edges are added to every $(i,j)$ unordered pair of discrete points, then the object is a simplex with $\binom {n} {2}$ 1-faces. If a subset of the $\binom {n} {2}$ possible edges is added, then special cases of emergent objects within the convex hull of ${\bf{X}}(0)$ may be obtained depending on how the edges are drawn. If no edges are added, then the complex object is array of points with coordinates that locate them in a $m$-dimensional space. The dimensions of this space may be geographical, or more generally dimensions of point-states on $m$ variables.  In the special case of a $1$-dimensional ${\bf{X}}(0)$ of point-states on one variable, the complex object is an array of points on the real number line, and this array may be displayed as a distribution of the number points located in positions (intervals) of the line.  I leave the application open and add no edges to the generated complex objects.  

Figures 1 and 2 are based on the same ${\bf{W}}$ and ${\bf{X}}(0)$. Figure 1 shows (a) a random array of initial positions (squares), (b) the simple object generated in the special case of ${\bf{A}} = {\bf{I}}$ (plus sign), and (c) the reduced compression effect of ${\bf{A}} = 0.80{\bf{I}}$ (solid circles). Figure 2 is based on a two-value $(0.10, 0.80)$ ${\bf{A}}$, which generates the reduced convergence of Figure 1 for a subset of $n$ and more reduced convergence for the remainder. The ${\bf{A}}$ values control the complex object within the convex hull of ${\bf{X}}(0)$, conditional on ${\bf{W}}$. 

Figures 3 and 4 display more orderly emergent complex objects that are resultants of particular combinations of ${\bf{A}}$ and ${\bf{X}}(0)$. Fewer initial positions than the size $n$ of the system appear when, as here, subsets of $n$ occupy the same initial positions.  

Figures 5 and 6 illustrate that if the initial space is organized as a polytope (here a triangle and a pentagon), with non-vacuous subsets of $n$ that occupy each vertex of the polytope, then the emergent complex object of the state-space process is constrained to the subspace of the polytope within the $m$-box.

Figure 7 displays a large scale $1$-dimensional application and its bar chart results. The result bears on Abelson's \cite{Abelson1964} unsolved problem in the field of opinion dynamics. His investigation of various models showed that a formal explanation of emergent consensus on specific issues is easily obtained, but that a formal explanation of emergent differentiated opinion clusters is a difficult, to his apparent consternation: ``Since universal ultimate agreement is an ubiquitous outcome of a very broad class of mathematical models, we are naturally led to inquire what on earth one must assume in order to generate the bimodal outcomes of community cleavage studies'' \cite[p.153]{Abelson1964}. One answer appears to be that such cleavage is based on an initial distribution of opinions on an issue concentrated around a moderate initial position that is then  ``pulled apart'' by a second order interpersonal influence process unfolding in a strong structure. Figure 8 displays the ${a_{ii}} \; i=1,...,n$ values that are involved in the process. The cleavage is based on a moderate mass of individuals who vary in their ${a_{ii}}$ values, and  extremists with  ${a_{ii}}$ values homogeneously near 1.

\section{Target Complex Objects}

Given $\{ {\bf{W}},{\bf{X}}(0), {\bf{X}}(\infty) \}$, each diagonal value of  ${\bf{A}}$, i.e., $a_{11},...,a_{nn}$, must satisfy equation \ref{eq:X0}. In scalar form, each $a_{ii}$ must simultaneously satisfy the constraints entailed in

\begin{equation}
0 < {a_{ii}} = \frac{{x_{ih}}(\infty)-{x_{ih}}(0)}
 {\sum_{j=1}^{n}{{w_{ij}} {x_{jh}}(\infty) - {x_{ih}}(0)}    } < 1, \;\; i=1,...,n; \;\; h=1,..,m.
 \label{eq:scalaraiisol}
\end{equation}

\noindent A specified ``target'' is feasible if and only if, for all $ i=1,...,n; h=1,..,m$, (a) the denominator is not zero, the numerator and denominators are of the same sign, and the absolute value of the numerator is less than or equal to the absolute value of the denominator. The domain of feasible ``targets'' is, therefore, subject to strong constraints, and it may be vacuous.

In contrast, given $\{ {\bf{W}}, {\bf{X}}(\infty) \}$, there exists an infinite number of  $\{{\bf{A}},{\bf{X}}(0)\}$ combinations with which the state-space process will converge to a specified complex object ${\bf{X}}(\infty)$. The infinite set of feasible $\{{\bf{A}} ,{\bf{X}}(0)\}$ combinations is defined by

\begin{equation}
{x_{ih}}(0) = \frac{{x_{ih}}(\infty) - {a_{ii}} \sum_{j=1}^{n}{ {w_{ij}} {x_{jh}}(\infty) }}{1-a_{ii}}, \;\; a_{ii} \in (0,1) .
\end{equation}

\noindent for $i=1,...,n; \;\; h=1,..,m$. With an arbitrary $a_{ii} \in (0,1)$ all values of ${x_{ih}}(0), \; h=1,...,m$ are determined for $i$. The set of arbitrary $a_{ii},\; i=1,...,n$ values fully specifies a feasible $\{{\bf{A}} ,{\bf{X}}(0)\}$ combination on the basis of which the state-space process will generate the specified  ${\bf{X}}(\infty)$. 

Whether particular classes of combinations have properties that are preferable to others will depend on the application. In this infinite set of feasible $\{{\bf{A}} ,{\bf{X}}(0)\}$ combinations, an attractive unbiased combination is obtained with  

\begin{equation}
{\bf{X}}(0) = ({\bf{I}} - {\bf{A}})^{-1}  \big( {\bf{I}} - {\bf{AW}} \big) {\bf{X}}(\infty) , \;\; {\bf{A}}  = (1/2){\bf{I}}.
\end{equation}  

\noindent This special case presents
 
\begin{equation}
{x_{ih}}(0) = 2{x_{ih}}(\infty) -  \sum_{j=1}^{n}{ w_{ij} {x_{jh}}(\infty)   },
\end{equation}
        
\noindent for which the values of ${\bf{X}}(0)$ are neither inflated or deflated by the choice of ${\bf{A}}$, and 

\begin{equation}
{x_{ih}}(\infty)  = \frac{ {x_{ih}}(0) -  \sum_{j=1}^{n}{ w_{ij} {x_{jh}}(\infty)   }}{2}
\end{equation}

\noindent for all $i=1,...,n; \;\; h=1,..,m$. However, heterogeneous $a_{11},...,a_{nn}$ values may be theoretically important.

An affine transformation of the specified system ${\bf{X}}(\infty) = {{\bf{V}}}{\bf{X}}(0)$ is permissible

\begin{equation}
\alpha + \beta {\bf{X}}(\infty) = {{\bf{V}}}(\alpha + \beta {\bf{X}}(0)),
\end{equation}

\noindent because the scalars $\{ \alpha, \beta \}$ pass through the system without altering ${\bf{V}}$.
 
\section{Discussion}

Analysis has been restricted to the case of aperiodic irreducible structures. This restriction may be relaxed to include the aperiodic reducible structures covered by a nonsingular ${\bf{I}} - {\bf{AW}}$, i.e, all aperiodic strong, unilateral, weak, and disconnected structures. 

In strong structures, the details of network typology do not affect the outcome; every node of a strong structure is a point basis of the network (a minimal subset of nodes from which all its nodes of the network are reachable). However, the point bases of unilateral and weak structures may be more complex, and the $a_{ii}$ dampening values of those nodes that constitute a point basis of the structure are especially important. In general, network typology affects the temporal trajectories of the state-space. Network typology also is an important factor in the sensitivity and vulnerability of complex objects to structural disturbances. Structural disturbances include changes of ${\bf{W}}$, ${\bf{A}}$, or ${\bf{X}}(0)$ during the state-space process.

Here, it may be noted that the $v_{ij}$ of ${\bf{V}}$ correspond to the relative net influence of node $j$ on node $i$. Thus, 

\begin{equation} 
{\bf{r}} = \frac{1}{n}{{\bf{V}}^{T}}{\bf{1}} , \;\; 
{{\bf{r}}^{T}} {\bf{1}} = 1, 
\end{equation}

\noindent is the average relative net influence of node $j$ on all nodes of the system, and may taken as a measure of the centrality of node $j$ in the system. In the special case of ${\bf{A}}={\bf{I}}$, where ${\bf{V}} = \mathop{\lim }\limits_{k \to \infty } {\bf{W}}^k$, ${\bf{r}}$ is the left eigenvector of ${\bf{W}}$ associated with $\lambda =  1$. In the special case of ${\bf{A}}=\alpha{\bf{I}}, \; \;  0< \alpha < 1$,  

\begin{align} 
\label{eq:TEC}
{\bf{r}} & = \frac{1}{n}{{\bf{V}}^{T}}{\bf{1}}, \\
& = \frac{1}{n}  ( {\bf{I}} - \alpha {\bf{W}}^{T})^{ - 1} {\bf{1}} ( 1 - \alpha ),  
 \notag
\end{align}

\noindent which may expressed as follows

\begin{equation} 
{\bf{r}}  = \frac{1-\alpha}{n} + \alpha{\bf{W}}^{T} {\bf{r}} 
\label{eq:PageRank}
\end{equation}

\noindent It is difficult for me to see how the PageRank equation \ref{eq:PageRank} may be treated as an invention when equation \ref{eq:TEC} was published prior to it \cite[p.1487]{Friedkin1991}.

I noted in the introduction that the convex hull of ${\bf{X}}(0)$ defines a space within which the initial coordinates of the $n$ points of the system include, as special cases, all convex polytopes in $m$ dimensions. I illustrated this in Figures 5 and 6, where the initial coordinates of the system presented a triangle and pentagon with each vertex occupied by a subset of $n$, and with all $i=1,...,n$ located in the vertices.  In such cases, all ${\bf{A}}$ will generate an end-state ${\bf{X}}(\infty)$ of coordinates, via the state-space process, that locates all $n$ points within the constrained subspace of the initial polytope. Conversely, it is interesting that solutions for ``target'' ${\bf{X}}(\infty)$ frequently present a reduced simplex for ${\bf{X}}(0)$ with fewer than $n$ vertices, i.e,  at least one vertex is occupied by multiple points with the same initial positions.

I have not addressed the temporal trajectory of ${\bf{X}}(k)$ as the state-space process alters the coordinates of the nodes. Effects of the detailed typology of ${\bf{W}}$ on these movements are of interest. Additionally, the structure of ${\bf{W}}$ (e.g., its point sets, cut sets, and disjoint path redundancies) is significant in determining the sensitivity and vulnerability of the system to in-process structural disturbances that alter the system's end-state. The present work suggests that in-process modifications of the nodal dampening values ${\bf{A}}$, which alter the values of the walks, may be used to readjust the system so that the points arrive at their specified destinations. 

Every walk of length $\ell$ in a network involves a path of lessor or equal length. As a particular $a_{ii}$ approaches zero, all paths involving edges from node $i$ are deactivated. In an irreducible network, the limits of a particular pair of values $({a_{ii}},{a_{jj}})$ present four binary combinations that control the existence of effects of the state of $i$ on the state $j$, and vice versa, for all $i$ and $j$. Generalizing to all $\binom {n} {r}$ combinations, ${\bf{A}}$ entirely controls the causal structure of direct and indirect influences of nodes' past states on other nodes' future states.           

A different implementation of the model (equation \ref{eq:FJ}) couples ${\bf{A}}$ and ${\bf{W}}$ with the assumption ${a_{ii}} = 1- {w_{ii}}$ for all $i$.  With such coupling, each node's resistance value ${w_{ii}}$, i.e., the structural value of its loop, corresponds to the extent to which the node is open or closed with respect to flows from it. Under this assumption ${\bf{A}} = {\bf{I}}$ corresponds to a ${\bf{W}}$ with a zero main diagonal, and ${\bf{A}} = {\bf{0}}$ corresponds to a ${\bf{W}}$ with a main diagonal of ones. Changes of ${\bf{A}}$ induce limited changes of ${\bf{W}}$.  The analysis that has been presented in this paper is fundamentally altered by this assumption, i.e., equation \ref{eq:scalaraiisol} is altered.  However, the associated analysis is tractable and reinforces the importance of the ${\bf{A}}$ construct in controlling system outcomes. I have not addressed this implementation in the present paper.

\bibliography{ArticleReferences}

 \begin{figure}[p]
 \centering
 \caption{One-value (0.80) A-matrix}
 \includegraphics[trim = 30mm 140mm 30mm 20mm, clip,width=0.8\textwidth]{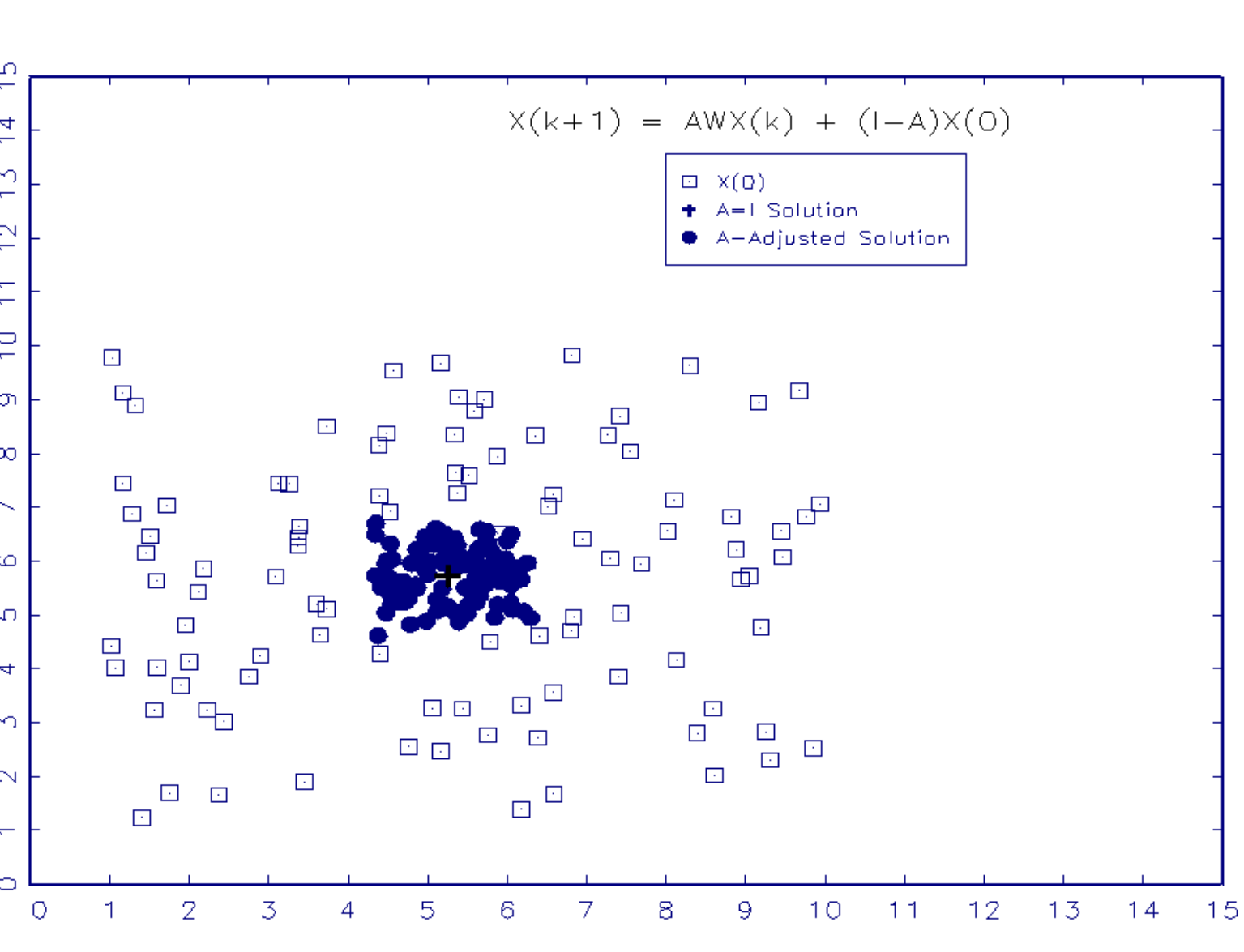}
 \caption{Two-value (0.10,0.80) A-matrix }
 \includegraphics[trim = 30mm 10mm 30mm 22mm, clip,width=0.8\textwidth]{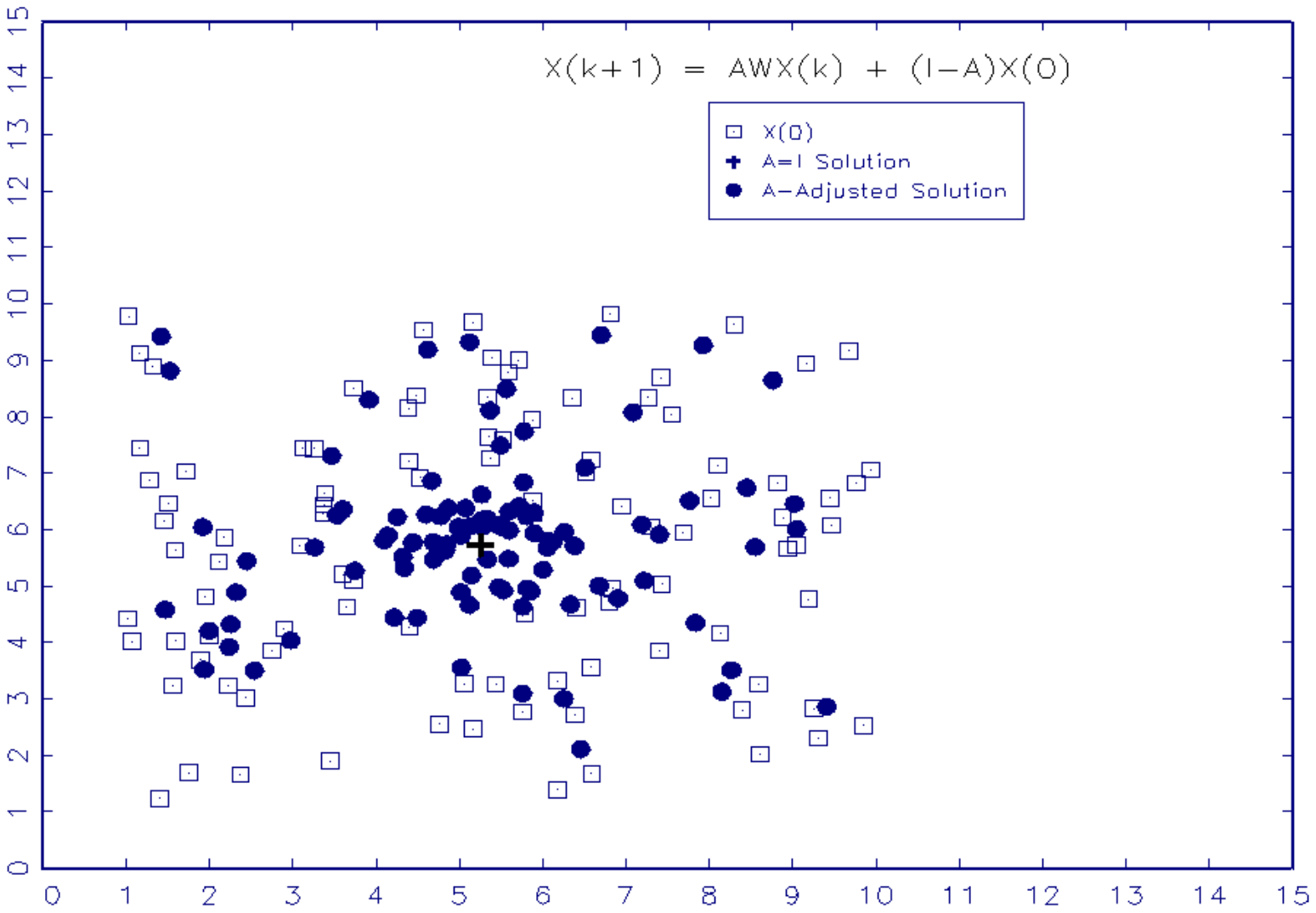}
 \end{figure}

 \begin{figure}[p]
 \centering
 \caption{More Orderly Complex Objects}
 \includegraphics[trim = 20mm 120mm 20mm 20mm, clip,width=0.8\textwidth]{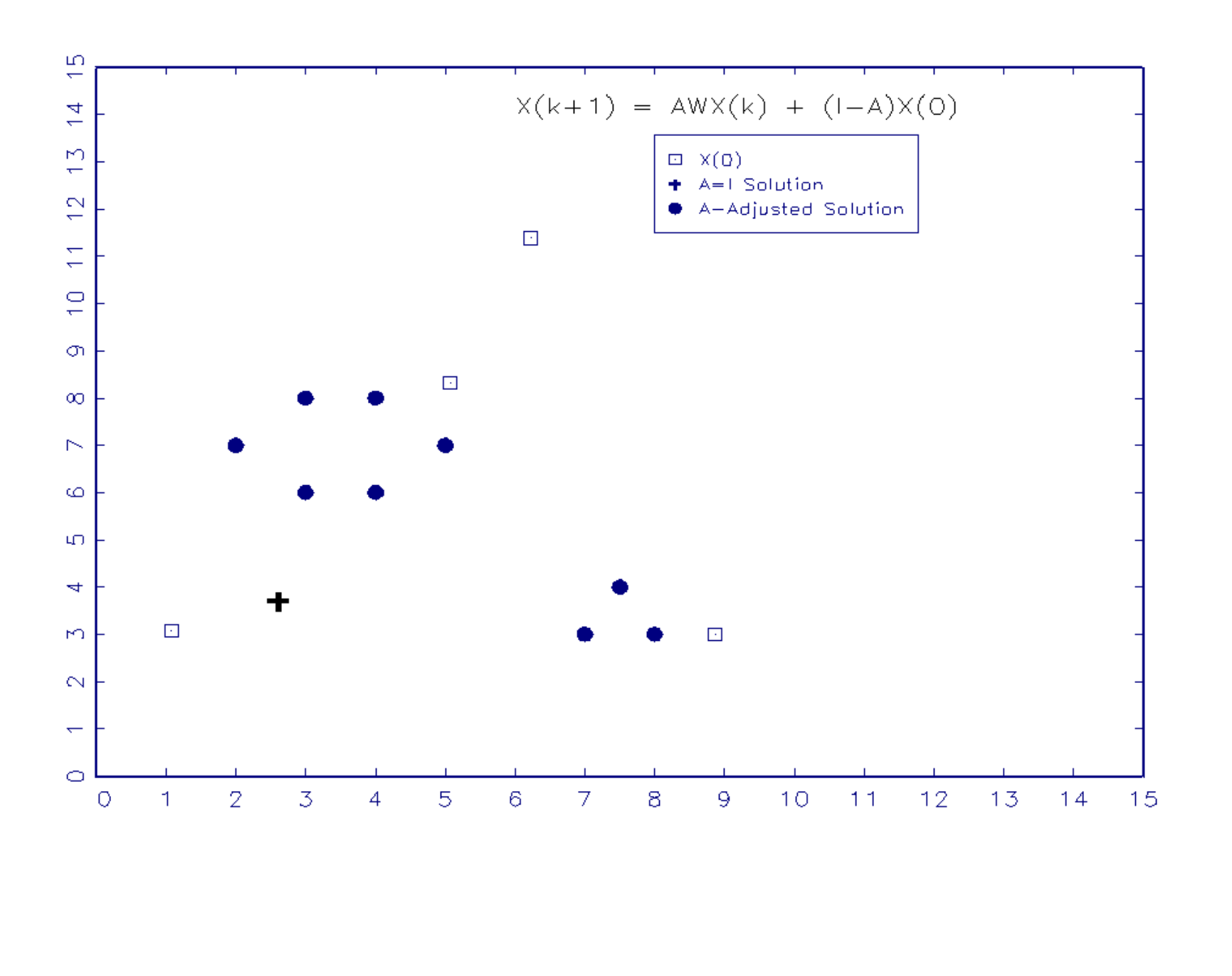}
 \caption{More Orderly Complex Objects }
 \includegraphics[trim = 29mm 20mm 10mm 20mm, clip,width=0.8\textwidth]{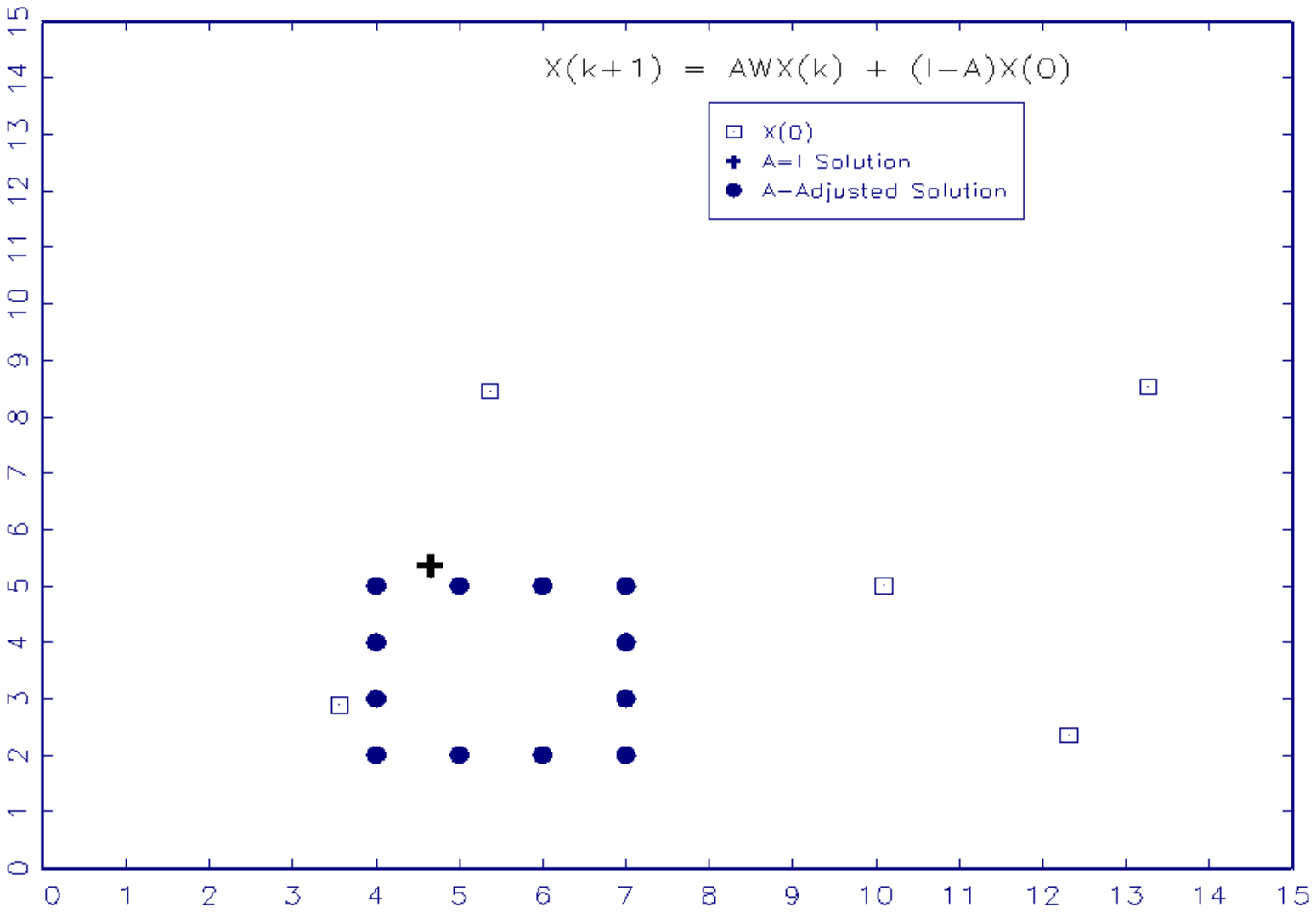}
 \end{figure}
 
 \begin{figure}[p]
  \centering
  \caption{Triangle Polytope Subspace Constraint}
  \includegraphics[trim = 20mm 120mm 20mm 20mm, clip,width=0.8\textwidth]{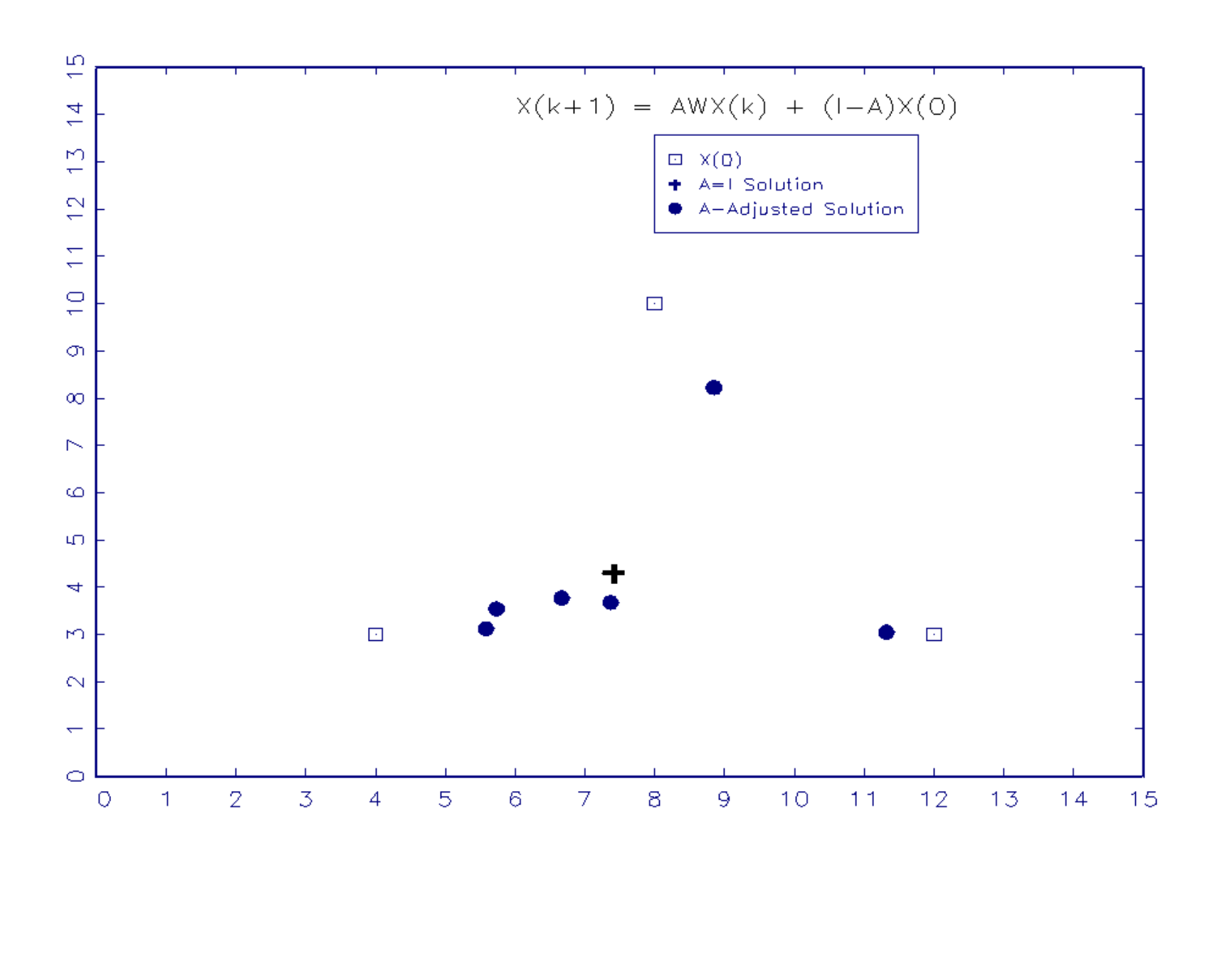}
  \caption{Pentagon Polytope Subspace Constraint}
  \includegraphics[trim = 20mm 20mm 10mm 20mm, clip,width=0.8\textwidth]{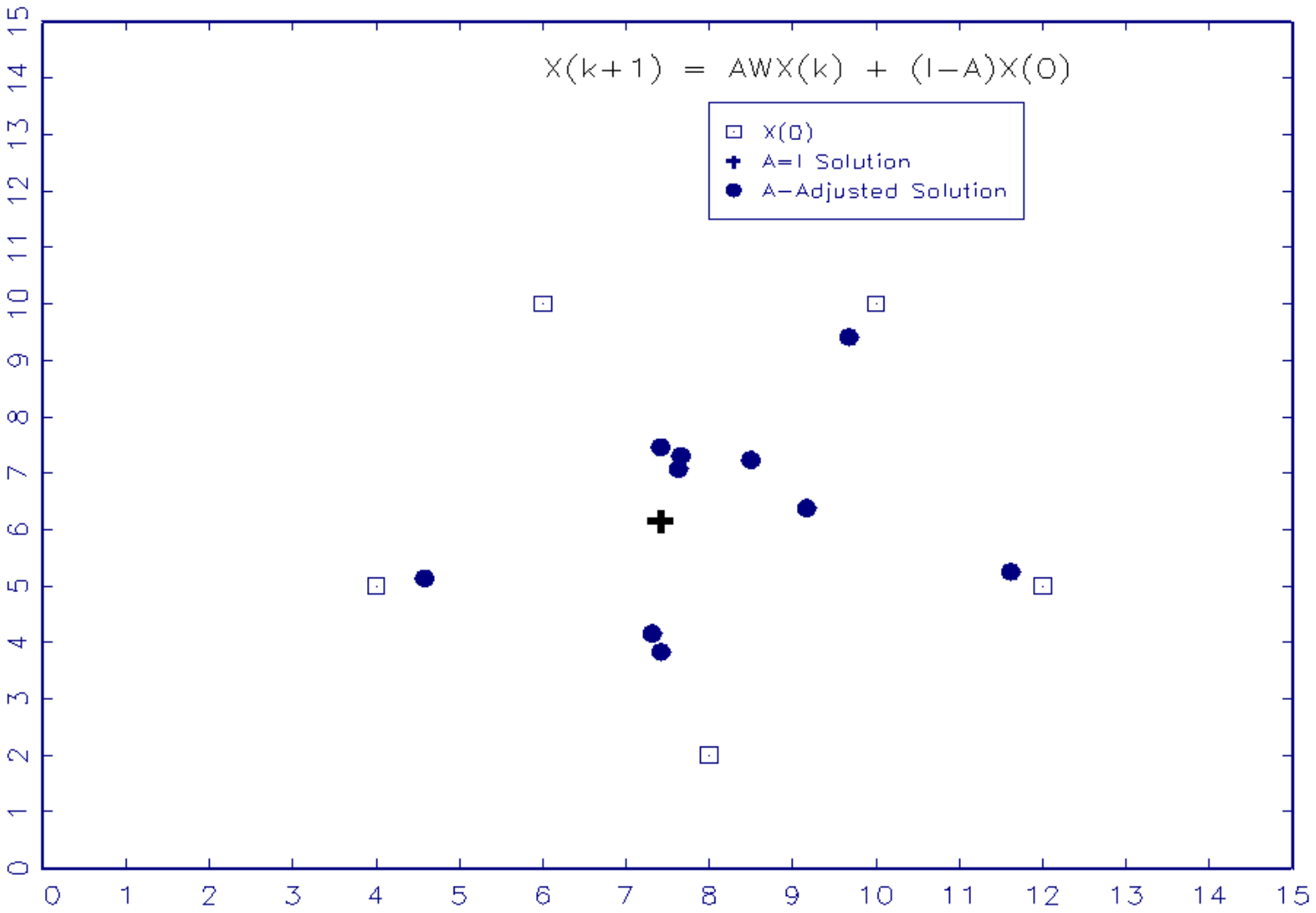}
  \end{figure}

 \begin{figure}[p]
 \centering
 \caption{Bimodal Complex Object \\ 
  (min$[{\bf{X}}(0)] = -10.756$; max$[{\bf{X}}(0)] = 13.655$; 
  min$[{\bf{X}}(\infty)] = 2.454$; max$[{\bf{X}}(\infty)] = 12.974$)
 }
 \includegraphics[trim = 20mm 120mm 20mm 20mm, clip,width=0.8\textwidth]{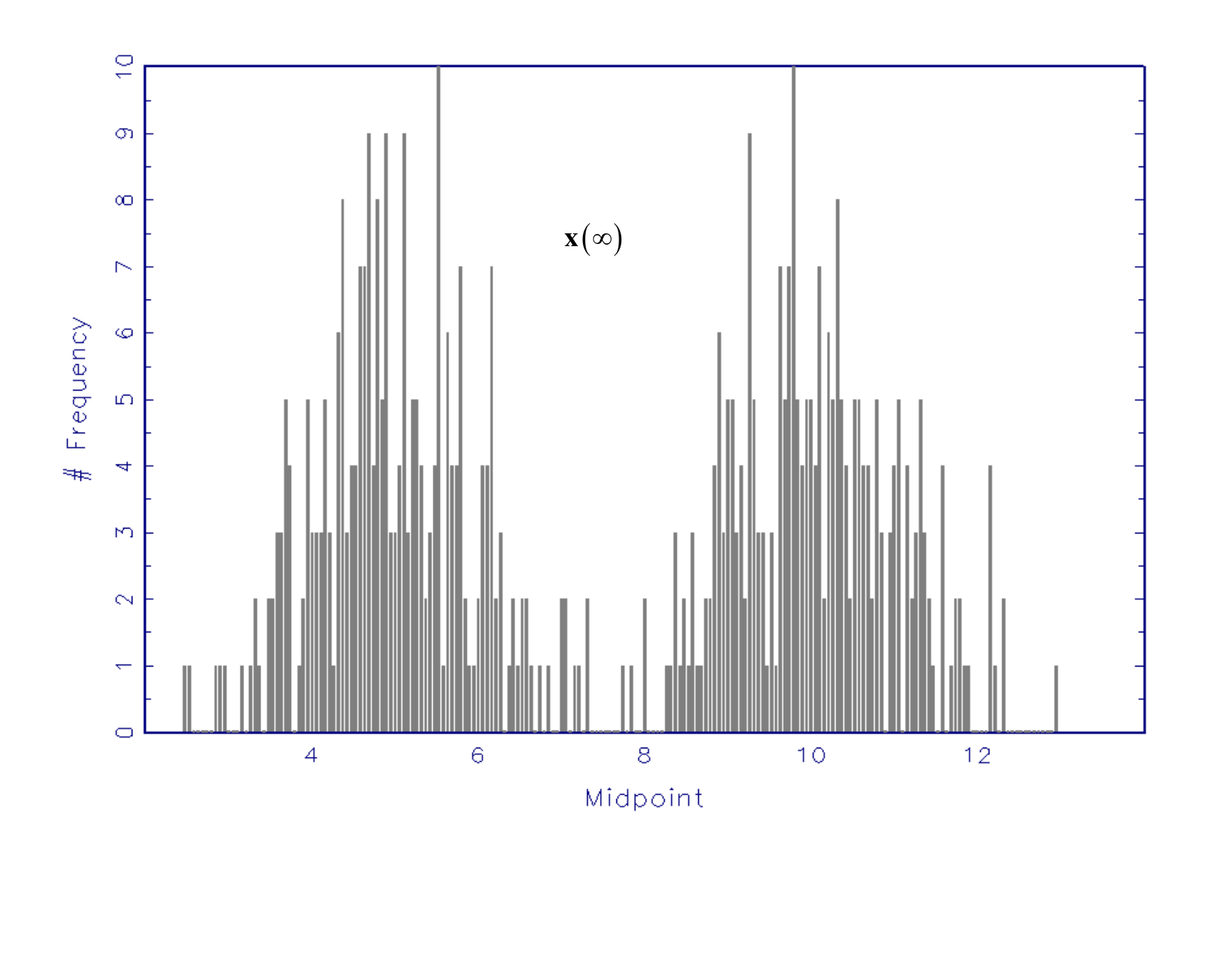}
  \includegraphics[trim = 20mm 20mm 10mm 20mm, clip,width=0.8\textwidth]{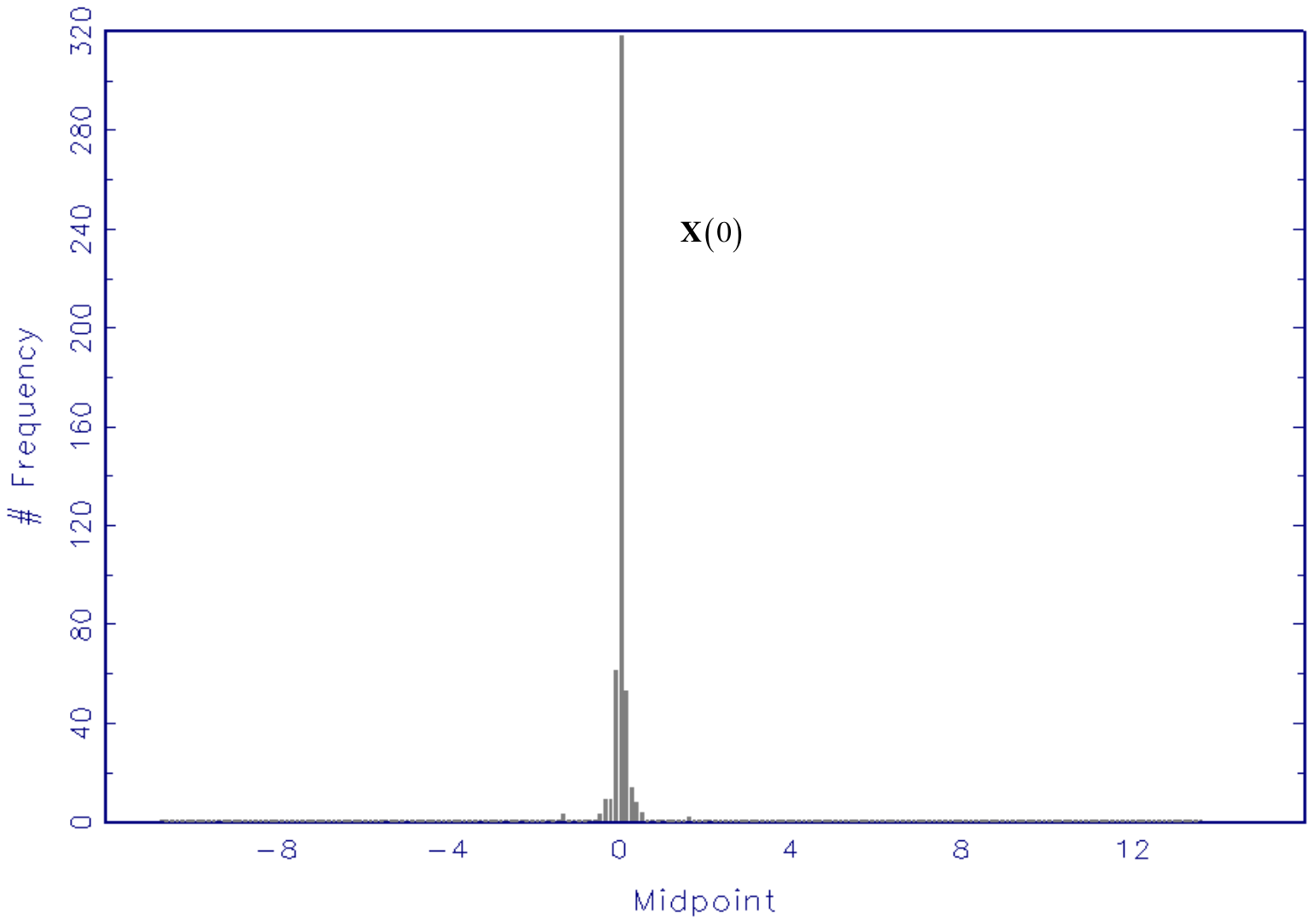}
  \label{fig:Bimodal}
  \end{figure}
 
 \begin{figure}[p]
  \centering
  \caption{Bimodal Complex Object}
  \includegraphics[trim = 20mm 20mm 10mm 20mm, clip,width=0.8\textwidth]{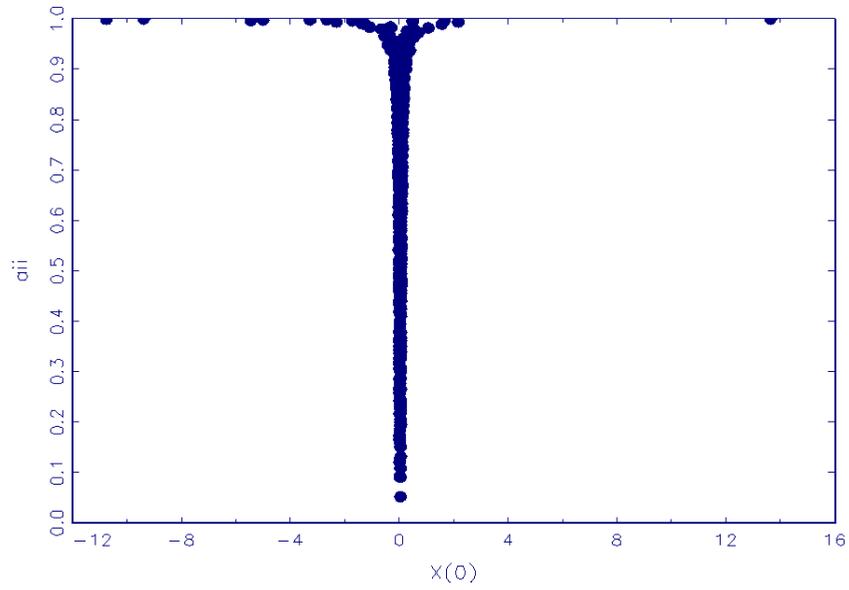}
  \label{fig:Bimodal}
  \end{figure}
  

\end{document}